\documentclass{article}

\usepackage{color}
\usepackage{amsmath}
\usepackage{amsfonts}
\usepackage{amssymb}
\usepackage{cancel}
\usepackage{authblk}

\newtheorem{lemma}{Lemma}

\begin{document}

\title{Cumulative probability for the sum of exponentially-distributed variables}
\author[1]{Cecilia Chirenti\footnote{email: cecilia.chirenti@ufabc.edu.br}}
\author[2]{M. Coleman Miller\footnote{email: miller@astro.umd.edu}}
\affil[1]{Centro de Matem\'{a}tica, Computa\c{c}\~{a}o e Cogni\c{c}\~{a}o, UFABC, 09210-170 Santo Andr\'{e}-SP, Brazil}
\affil[2]{Department of Astronomy and Joint Space-Science Institute, University of Maryland, College Park, MD 20742-2421 USA}
\date{}

\maketitle

\begin{abstract}

Exponential distributions appear in a wide range of applications including chemistry, nuclear physics, time series analyses, and stock market trends.  There are conceivable circumstances in which one would be interested in the cumulative probability distribution of the sum of some number of exponential variables, with potentially differing constants in their exponents.  In this article we present a pedagogical derivation of the cumulative distribution, which reproduces the known formula from power density analyses in the limit that all of the constants are equal, and which assumes no prior knowledge of combinatorics except for some of the properties of a class of symmetric polynomials in $n$ variables (Schur polynomials).

\end{abstract}


\section{Introduction}

In many circumstances, a system can be considered to have no ``memory", in the sense that, for example, the expected amount of time to the next radioactive decay, or chemical reaction, or failure of a machine part, does not depend on how much time has already elapsed.  The unique memoryless and continuous probability distribution is the exponential distribution; if a variable $P$ has an expected value $P_{\rm exp}$, then the normalized differential probability density is
\begin{equation}
{\rm prob}(P)={1\over{P_{\rm exp}}}e^{-P/P_{\rm exp}}
\end{equation}
and the probability that $P$ will be $P_0$ or larger is
\begin{equation}
{\rm prob}(P\geq P_0)=e^{-P_0/P_{\rm exp}}\; .
\end{equation}
Alternatively, the exponential distribution is the solution to the differential equation
\begin{equation}
{dy\over{dP}}=-y/P_{\rm exp}\; .
\end{equation}

There may be occasions in which one wishes to compute the probability distribution of the sum of several independently-sampled exponential variables, which could have different decay constants.  The sum of the waiting times for a decay from several different radioactive substances provides a perhaps artificial example of this type.  

Thus we are faced with the following problem: if we combine the {\it statistically independent} measurements of a quantity $P$ for $n$ discrete variables, what is the probability distribution of the sum of the values of those quantities given that they are all exponentially distributed?  Put another way, if the expected constants for the $n$ variables are $P_{1,{\rm exp}}, P_{2,{\rm exp}},\ldots,P_{n,{\rm exp}}$, how probable is it in that model that the sum of the measured values will be at least $P_{\rm tot}$?  

Motivated by this problem, we will proceed as follows. In Section \ref{sec:integral} we obtain by mathematical induction our general formula. This proof is completed by a basic lemma of long standing, presented in Section \ref{sec:lemma}, in which we use some properties of Schur polynomials. In Section \ref{sec:application} we apply our formula to the case where all exponential constants are the same, and recover the correct result, which is known in power density theory among other contexts. We present our final remarks in Section \ref{sec:final}.

\section{The general formula}
\label{sec:integral}

The probability that the sum of the values of $n$ measurements drawn from exponentials is greater than or equal to $P_{\rm tot}$, is one minus the probability that the sum is less than $P_{\rm tot}$: 
\begin{equation}
\begin{array}{rl}
   {\rm prob}(P \geq P_{\rm tot})&=1-\int_0^{P_{\rm tot}} (1/P_{1,{\rm exp}})e^{-P_1/P_{1,{\rm exp}}} \\
                   &\times \biggl[\int_0^{P_{\rm tot}-P_1} (1/P_{2,{\rm exp}})e^{-P_2/P_{2,{\rm exp}}} \\
                   &\times \biggl[\int_0^{P_{\rm tot}-P_1-P_2} (1/P_{3,{\rm exp}})e^{-P_3/P_{3,{\rm exp}}}\times\ldots\times \\                   
                   &\times\biggl[\int_0^{P_{tot}-\sum_{i=1}^{n-1}P_i} (1/P_{n,{\rm exp}})e^{-P_n/P_{n,{\rm exp}}}dP_n\biggr]dP_{n-1}\biggr]\ldots\biggr]dP_1\; .
\end{array}
\label{eq:original}
\end{equation}
We claim that this expression evaluates to
\begin{equation}
{\rm prob}(P \geq P_{\rm tot})=\sum_{i=1}^n \left({P_{i,{\rm exp}}^{n-1}e^{-P_{\rm tot}/P_{i,{\rm exp}}}\over{\prod_{j\neq i}(P_{i,{\rm exp}}-P_{j,{\rm exp}})}}\right)\; .
\label{eq:formula}
\end{equation}

We prove this by induction:

{\bf Step 1:} For $n=1$, the probability is
\begin{equation}
\begin{array}{rl}
{\rm prob}(P \geq P_{\rm tot})&=1-\int_0^{P_{\rm tot}}(1/P_{1,{\rm exp}})e^{-P_1/P_{1,{\rm exp}}}dP_1\\
&=1-(1-e^{-P_{\rm tot}/P_{1,{\rm exp}}})=e^{-P_{\rm tot}/P_{1,{\rm exp}}}\; .
\end{array}
\end{equation}
This is consistent with eq. (\ref{eq:formula}).

{\bf Step 2:} Assume that the formula is true for a given $n$. For $n+1$, insertion of the formula along with some relabeling gives us
\begin{eqnarray}
&&{\rm prob}(P \geq P_{\rm tot})=1-\int_0^{P_{\rm tot}}{1\over{P_{1,{\rm exp}}}}e^{-P_1/P_{1,{\rm exp}}}dP_1\times\nonumber\\
&\times&\left[1-\sum_{i=2}^{n+1}{P_{i,{\rm exp}}^{n-1}e^{-(P_{\rm tot}-P_1)/P_{i,{\rm exp}}}\over{\prod_{1<j\neq i}(P_{i,{\rm exp}}-P_{j,{\rm exp}})}}\right]\; .
\end{eqnarray}
Multiplying through, we obtain
\begin{eqnarray}
&&1-\left[\int_0^{P_{\rm tot}}{1\over{P_{1,{\rm exp}}}}e^{-P_1/P_{1,{\rm exp}}}dP_1-\right.\nonumber\\
&-&\left.\int_0^{P_{\rm tot}}\sum_{i=2}^{n+1}{e^{-P_1(1/P_{1,{\rm exp}}-1/P_{i,{\rm exp}})}P_{i,{\rm exp}}^{n-1}e^{-P_{\rm tot}/P_{i,{\rm exp}}}\over{P_{1,{\rm exp}}\prod_{1<j\neq i}(P_{i,{\rm exp}}-P_{j,{\rm exp}})}}dP_1\right]\; .
\end{eqnarray}
The first integral just evaluates to $1-e^{-P_{\rm tot}/P_{1,{\rm exp}}}$.  For the second integral, we note that the factor in the exponent multiplying $-P_1$ is
\begin{equation}
{1\over{P_{1,{\rm exp}}}}-{1\over{P_{i,{\rm exp}}}}={P_{i,{\rm exp}}-P_{1,{\rm exp}}\over{P_{1,{\rm exp}}P_{i,{\rm exp}}}}
\end{equation}
and thus the second integral becomes
\begin{eqnarray}
&-&\sum_{i=2}^{n+1}{P_{i,{\rm exp}}^ne^{-P_{\rm tot}/P_{i,{\rm exp}}}\over{P_{i,{\rm exp}}-P_{1,{\rm exp}}}}{1\over{\prod_{1<j\neq i}(P_{i,{\rm exp}}-P_{j,{\rm exp}})}}\times\nonumber\\
&\times&\left[1-e^{-P_{\rm tot}/P_{1,{\rm exp}}}e^{P_{\rm tot}/P_{i,{\rm exp}}}\right]\; .
\end{eqnarray}
Subtracting the sum of the integrals from 1 gives
\begin{eqnarray}
&&{\rm prob}(P \geq P_{\rm tot})=\sum_{i=2}^{n+1}{P_{i,{\rm exp}}^ne^{-P_{\rm tot}/P_{i,{\rm exp}}}\over{\prod_{j\neq i}(P_{i,{\rm exp}}-P_{j,{\rm exp}})}}+e^{-P_{\rm tot}/P_{1,{\rm exp}}}\times\nonumber\\
&\times&\left[1-\sum_{i=2}^{n+1}{P_{i,{\rm exp}}^n\over{\prod_{j\neq i}(P_{i,{\rm exp}}-P_{j,{\rm exp}})}}\right]\; .
\end{eqnarray}
The first term has the required form for $i=2$ through $i=n+1$.  We will prove our formula if we can show that
\begin{equation}
1-\sum_{i=2}^{n+1}{P_{i,{\rm exp}}^n\over{\prod_{j\neq i}(P_{i,{\rm exp}}-P_{j,{\rm exp}})}}={P_{1,{\rm exp}}^n\over{\prod_{j\neq 1}(P_{1,{\rm exp}}-P_{j,{\rm exp}})}}\; ,
\end{equation}
which is to say that we can prove our formula if we can show that
\begin{equation}
\sum_{i=1}^{n+1}{P_{i,{\rm exp}}^n\over{\prod_{j\neq i}(P_{i,{\rm exp}}-P_{j,{\rm exp}})}}=1\; .
\label{eq:lemma}
\end{equation}

\section{Statement and proof of the basic lemma}
\label{sec:lemma}

Now we proceed to prove eq. (\ref{eq:lemma}) in order to complete the proof by induction of our general formula given in eq. (\ref{eq:formula}). As this result will be used again in Section \ref{sec:application}, we present it in this Section as a separate lemma (which, as we discuss, already appeared in generalized form some centuries ago), which we state and prove below using techniques that do not require advanced mathematics. We give explicit examples with $n = 1,2$ and 3 in Appendix \ref{sec:ex-lemma}.

\begin{lemma}
\label{lemma1}
Let $\{x_i\}_{i \in \mathbb{N}}$ be a set of real positive numbers, $n$ a given natural number and $\{a_i\}_{i \in \mathbb{N}}$ be a real sequence given by
\begin{equation}
a_i = \frac{x_i^n}{\prod_{j\neq i}(x_i - x_j)}\,, \quad \rm{with} \quad j = 1,\ldots,i-1,i+1,\ldots,n+1\,.
\end{equation}
Then the sum $S_{n+1} = \sum_{i=1}^{n+1} a_i = 1$.
\end{lemma}
As we indicate above, generalizations of this formula have been applied to interpolation theory as early as \cite{Newton}, and there are significant connections to geometry (in evaluating integrals over projection space and in the context of flag manifolds; for several papers that clarify these connections see \cite{LS1,LS2,LS3,LS4}, and for a more recent review see \cite{Pragacz}).  In addition, when one interprets the exponential as a probability distribution, one can take advantage of the fact that the distribution of the sum of two probability densities is the convolution of the densities (for an application in this context, see, e.g., \cite{Akkouchi}).  Here, however, we present a proof that relies on the properties of certain classes of symmetric polynomials.

{\it Proof.} $S_{n+1}$ is a symmetric polynomial in the $x_i$ (in fact, it is a Schur polynomial as we will see below) because it can be written as $A_n/B_n$, where both $A_n$ and $B_n$ are alternating polynomials of degree $n$  in the $x_i$'s and are given by
\begin{eqnarray}
\label{eq:A_n}
A_n &=& \sum_{i=1}^{n+1}(-1)^{i+1}x_i^n\prod_{j < k; j,k \neq i}(x_j-x_k)\,,\\
B_n &=& \prod_{j < k}(x_j-x_k)\,.
\label{eq:B_n}
\end{eqnarray}
%
All alternating polynomials can be written as the product $V_n\cdot s$, where $V_n$ is the Vandermonde polynomial (\textit{i.e.}, the determinant of the Vandermonde matrix \cite{NumericalRecipes}) and $s$ is a symmetric polynomial. We recognize now that $B_n = V_n$, and therefore $S_{n+1}$ is a symmetric polynomial given by the ratio of the alternating polynomial $A_n$ over the Vandermonde polynomial. Thus, $S_{n+1}$ is called a Schur polynomial \cite{Schur}.

Therefore we can write $A_n = B_n\cdot s$. As $A_n$ and $B_n$ are polynomials of the same order, $s$ must be a zeroth-order polynomial: $s = C$, with $C = $ const. and $A_n = C\cdot B_n$. To show that $C = 1$, we compare the coefficients of the $x_1^n$ terms in $A_n$ and $B_n$, as given by eqs. (\ref{eq:A_n}) and (\ref{eq:B_n}). For both $A_n$ and $B_n$, the coefficient of $x_1^n$ is
\begin{equation}
\prod_{1< j < k}(x_j-x_k)\,,
\end{equation}
which demonstrates that $C = 1$. Therefore, $A_n = B_n$ and $S_{n+1} = A_n/B_n = 1$. 
 
\section{An application: Reduction of the general formula in the case that all expected values are nearly equal}
\label{sec:application}

Our general formula given by eq. (\ref{eq:formula}) has, in its denominator, the products of differences of expected values.  If two or more of the $P_{i,{\rm exp}}$ approach each other, then the formula appears singular.  It is not, in fact, as the original integral given in eq. (\ref{eq:original}) is clearly nonsingular. We will demonstrate this below, and we give explicit examples for $n = 2,3$ in Appendix \ref{sec:ex-app}.

Let us assume that there are $n$ total measurements, and that all of them have expected values that are very close to each other. We will assume that $P_{2,{\rm exp}}=P_{1,{\rm exp}}+\epsilon_2$, $P_{3,{\rm exp}}=P_{1,{\rm exp}}+\epsilon_3$, $\ldots$, $P_{n,{\rm exp}}=P_{1,{\rm exp}}+\epsilon_n$, and take the limit of the resulting expression as $\epsilon_2\rightarrow 0$, $\epsilon_3\rightarrow 0,\ldots, \epsilon_m\rightarrow 0$:
\begin{eqnarray}
&&\lim_{P_{i,{\rm exp}}\to P_{1,{\rm exp}}}{\rm prob}(P \geq P_{\rm tot}) = \nonumber\\
&=&\lim_{\epsilon_i \to 0}\left(\frac{P_{1,{\rm exp}}^{n-1}e^{-P_{\rm tot}/P_{1,{\rm exp}}}}{\prod_{j>1}(-1)^{n-1}\epsilon_j} 
+ \sum_{i=2}^n\frac{(P_{1,{\rm exp}}+\epsilon_i)^{n-1}e^{-P_{\rm tot}/(P_{1,{\rm exp}}+\epsilon_i)}}{\epsilon_i\prod_{j > 1; j \neq i}(\epsilon_i - \epsilon_j)}\right)
\end{eqnarray}
or
\begin{equation}
\lim_{\epsilon_i \to 0}\left(\frac{f(P_{1,{\rm exp}})}{\prod_{j>1}(-1)^{n-1}\epsilon_j} 
+ \sum_{i=2}^n\frac{f(P_{1,{\rm exp}}+\epsilon_i)}{\epsilon_i\prod_{j > 1; j \neq i}(\epsilon_i - \epsilon_j)}\right)\,, 
\label{eq:lim_f}
\end{equation}
where we used the auxiliary function
\begin{equation}
f(P) = P^{n-1}e^{-P_{\rm tot}/P}\,.
\label{eq:f}
\end{equation}
We will calculate this limit by expanding $f(P_{1,{\rm exp}}+\epsilon_i)$ in a Taylor series around $P_{1,{\rm exp}}$,
\begin{equation}
f(P_{1,{\rm exp}}+\epsilon_i) = \sum_{l=0}^{\infty}\frac{f^{(l)}(P_{1,{\rm exp}})}{l!}\epsilon_i^l\,,
\end{equation}
and substituting it in eq. (\ref{eq:lim_f}). We can now collect the terms proportional to $f(P_{1,{\rm exp}})$ and its  $l$-derivatives.

For $l=0$, we find the coefficient of $f(P_{1,{\rm exp}})$ and we can show that it is:
\begin{equation}
\frac{1}{\prod_{j>1}(-1)^{n-1}\epsilon_j} 
+ \sum_{i=2}^n\frac{1}{\epsilon_i\prod_{j > 1;j \neq i}(\epsilon_i - \epsilon_j)} = 0\,, 
\label{eq:k0}
\end{equation}
where, similarly to our proof of Lemma \ref{lemma1}, we can write the second term as
\begin{equation}
\frac{\sum_{i=2}^n(-1)^{i}\prod_{j \neq i}\epsilon_j\prod_{1<j<k;j,k \neq i}(\epsilon_j-\epsilon_k)}{\prod_{i>1}\epsilon_i\prod_{1<j < k}(\epsilon_j - \epsilon_k)} = \frac{(-1)^n}{\prod_{i>1}\epsilon_i}\,, 
\end{equation}
after the cancellation of the alternating polynomials of degree $n-2$ in the $\epsilon_j$'s, leading to the result in (\ref{eq:k0}).

For each $l$ in the range $1< l < n-1$, we find that the coefficient of the l-derivative $f^{(l)}(P_{1,{\rm exp}})$ is:
\begin{equation}
\frac{1}{l!}\sum_{i=2}^n\frac{\epsilon_i^l}{\prod_{j \neq 1}\epsilon_i(\epsilon_i - \epsilon_j)} = 0\,, 
\end{equation}
because the sum can be written in the form 
\begin{equation}
\frac{\sum_{i=2}^n(-1)^{i}\epsilon_i^{l-1}\prod_{1<j<k;j,k \neq i}(\epsilon_j-\epsilon_k)}{\prod_{1<j < k}(\epsilon_j - \epsilon_k)}\; .
\end{equation}
From this we can see that the denominator is the Vandermonde determinant in the variables $\epsilon_i$ (with $i = 2, \ldots,n$) \cite{NumericalRecipes}, whereas the numerator is zero, because it is equal to the determinant of the modified Vandermonde matrix in which the line $\epsilon_i^{l-1}$ appears twice, where in the second occurrence it replaces the last line $\epsilon_i^{n-2}$.

For $l = n-1$, the coefficient of the l-derivative $f^{(l)}(P_{1,{\rm exp}})$ is
\begin{equation}
\frac{1}{(n-1)!}\sum_{i=2}^n\frac{\epsilon_i^{(n-1)}}{\prod_{j \neq 1}\epsilon_i(\epsilon_i - \epsilon_j)} = \frac{1}{(n-1)!}\,, 
\end{equation}
where we can again use Lemma \ref{lemma1} to show that the sum equals 1.

For $l > n-1$, all coefficients of the $l$-derivatives $f^{(l)}(P_{1,{\rm exp}})$ go to zero when we take the limit $\epsilon_i \to 0$, as these coefficients will have higher order polynomials in $\epsilon_i$ in the numerator.

Thus, finally, we find that the desired result for the limit given in eq. (\ref{eq:lim_f}) is given by the only surviving term  $f^{(n-1)}(P_{1,{\rm exp}})/(n-1)!$ and, for the derivatives of our $f(P)$ given by eq. (\ref{eq:f}) we prove in Appendix \ref{sec:derivative} that
\begin{equation}
f^{(n-1)}(P_{1,{\rm exp}}) = (n-1)!\, e^{-P_{\rm tot}/P_{1,{\rm exp}}}\sum_{i=0}^{n-1}\frac{(P_{\rm tot}/P_{1,{\rm exp}})^i}{i!}\,,
\label{eq:df}
\end{equation}
so that we obtain 
\begin{equation}
\lim_{P_{i,{\rm exp}}\to P_{1,{\rm exp}}}{\rm prob}(P \geq P_{\rm tot}) = e^{-P_{\rm tot}/P_{1,{\rm exp}}}\sum_{i=0}^{n-1}\frac{(P_{\rm tot}/P_{1,{\rm exp}})^i}{i!}\,.
\end{equation}
This is the result known from, for example, time series analysis\footnote{If we apply Leahy normalization, as is standard in astronomy, then $P_{\rm exp}=2$ and thus the probability is
$e^{-P_{\rm tot}/2}\sum_{i=0}^{n-1}{(P_{\rm tot}/2)^i\over{i!}}\,$.} in the context of an intrinsically steady source whose only power comes from Poisson fluctuations (see \cite{BlackmanTukey} and the limit $P_s\rightarrow 0$ of equation (16) from \cite{Groth75}).

\section{Final Remarks}
\label{sec:final}

We have presented here an application of basic properties of Schur polynomials (symmetric polynomials that can be written as the ratio between an alternating polynomial in the numerator and the Vandermonde polynomial in the denominator) motivated by the statistical problem of determining the cumulative distribution of the sum of exponentially-distributed variables. Our work led to the result given in Lemma \ref{lemma1}, for calculating the sums of  Schur polynomials that appeared in our analysis.

As an application of our analysis, we recovered the known expression for the probability that the power density of an intrinsically constant source, summed over some number of independent frequency bins, exceeds a given total. As part of our calculation, we found an interesting direct expression for the $(n-1)$th-derivative of our auxiliary function $f$.

\medskip
{\bf Acknowledgments}
We are very grateful to our colleagues Pedro Lauridsen Ribeiro and Harry Tamvakis, who read through an earlier version of this paper and provided invaluable context.  This work was supported in part by S{\~{a}}o Paulo Research Foundation (FAPESP) grant 2015/20433-4 and by joint research workshop award 2015/50421-8 from FAPESP and the University of Maryland, College Park.

\appendix

\section{Examples of expressions using the basic lemma}
\label{sec:ex-lemma}

Here we present some examples of application of Lemma \ref{lemma1} for different values of $n$.

For $n = 1$:
\begin{equation}
\frac{x_1}{x_1-x_2} + \frac{x_2}{x_2-x_1} = \frac{x_1-x_2}{x_1-x_2} = 1\,.
\label{eq:lemma1-n1}
\end{equation}

For $n = 2$:
\begin{eqnarray}
&&\frac{x_1^2}{(x_1-x_2)(x_1-x_3)} + \frac{x_2^2}{(x_2-x_1)(x_2-x_3)} + \frac{x_3^2}{(x_3-x_1)(x_3-x_2)} = \nonumber\\
&=& \frac{x_1^2(x_2-x_3)-x_2^2(x_1-x_3)+x_3^2(x_1-x_2)}{(x_1-x_2)(x_1-x_3)(x_2-x_3)} = 1\,.
\end{eqnarray}

For $n = 3$:
\begin{eqnarray}
&&\frac{x_1^3}{(x_1-x_2)(x_1-x_3)(x_1-x_4)} + \frac{x_2^3}{(x_2-x_1)(x_2-x_3)(x_2-x_4)} + \nonumber\\
&+&\frac{x_3^3}{(x_3-x_1)(x_3-x_2)(x_3-x_4)} + \frac{x_4^3}{(x_4-x_1)(x_4-x_2)(x_4-x_3)} = \nonumber\\
&=& \frac{x_1^3(x_2-x_3)(x_2-x_4)(x_3-x_4)-x_2^3(x_1-x_3)(x_1-x_4)(x_3-x_4)}{(x_1-x_2)(x_1-x_3)(x_1-x_4)(x_2-x_3)(x_2-x_4)(x_3-x_4)} + \nonumber\\
&+& \frac{x_3^3(x_1-x_2)(x_1-x_4)(x_2-x_4)-x_4^3(x_1-x_2)(x_1-x_3)(x_2-x_3)}{(x_1-x_2)(x_1-x_3)(x_1-x_4)(x_2-x_3)(x_2-x_4)(x_3-x_4)} \nonumber\\
&=& 1\,.
\end{eqnarray}

\section{Examples of the reduction of the general formula}
\label{sec:ex-app}

Here we present some examples of the reduction of the general formula (\ref{eq:formula}) in the case that all $n$ power densities are equal, for different values of $n$.

For $n = 2$ we have the simplest case
\begin{eqnarray}
&&\lim_{P_{2,{\rm exp}}\to P_{1,{\rm exp}}}\left(\frac{P_{1,{\rm exp}}e^{-P_{\rm tot}/P_{1,{\rm exp}}}}{P_{1,{\rm exp}}-P_{2,{\rm exp}}} +
\frac{P_{2,{\rm exp}}e^{-P_{\rm tot}/P_{2,{\rm exp}}}}{P_{2,{\rm exp}}-P_{1,{\rm exp}}}\right) = \nonumber\\
&=& \lim_{\epsilon_2 \to 0}\left(\frac{P_{1,{\rm exp}}e^{-P_{\rm tot}/P_{1,{\rm exp}}}}{-\epsilon_2} + 
\frac{(P_{1,{\rm exp}}+\epsilon_2)e^{-P_{\rm tot}/P_{1,{\rm exp}}}}{\epsilon_2}\right) \nonumber\\
&=&  \lim_{\epsilon_2 \to 0} \frac{f(P_{1,{\rm exp}}+\epsilon_2)-f(P_{1,{\rm exp}})}{\epsilon_2} \nonumber\\
&=& f'(P_{1,{\rm exp}}) 
= e^{-P_{\rm tot}/P_{1,{\rm exp}}}\left(1 + \frac{P_{\rm tot}}{P_{1,{\rm exp}}} \right)\,, 
\end{eqnarray}
where we used $P_{2,{\rm exp}} = P_{1,{\rm exp}} + \epsilon_2$ and $f(P) = Pe^{-P_{\rm tot}/P}$.

For $n = 3$, we have
\begin{eqnarray}
&&
\lim_{\begin{array}{c}
{\scriptstyle P_{2,{\rm exp}}\to P_{1,{\rm exp}}}\\
{\scriptstyle P_{3,{\rm exp}}\to P_{1,{\rm exp}}}
\end{array}}
\left(\frac{P_{1,{\rm exp}}^2e^{-P_{\rm tot}/P_{1,{\rm exp}}}}{(P_{1,{\rm exp}}-P_{2,{\rm exp}})(P_{1,{\rm exp}}-P_{3,{\rm exp}})} 
+\right. \nonumber\\
&+&\left. 
\frac{P_{2,{\rm exp}}^2e^{-P_{\rm tot}/P_{2,{\rm exp}}}}{(P_{2,{\rm exp}}-P_{1,{\rm exp}})(P_{2,{\rm exp}}-P_{3,{\rm exp}})} +
\frac{P_{3,{\rm exp}}^2e^{-P_{\rm tot}/P_{2,{\rm exp}}}}{(P_{3,{\rm exp}}-P_{1,{\rm exp}})(P_{3,{\rm exp}}-P_{2,{\rm exp}})}\right) = \nonumber\\
&=& \lim_{\epsilon_2,\epsilon_3\to 0}
\left(\frac{P_{1,{\rm exp}}^2e^{-P_{\rm tot}/P_{1,{\rm exp}}}}{\epsilon_2\epsilon_3} + 
\frac{(P_{1,{\rm exp}}+\epsilon_2)^2e^{-P_{\rm tot}/P_{1,{\rm exp}}}}{\epsilon_2(\epsilon_2-\epsilon_3)} +\right.
\nonumber\\
&+&\left.\frac{(P_{1,{\rm exp}}+\epsilon_3)^2e^{-P_{\rm tot}/P_{1,{\rm exp}}}}{\epsilon_3(\epsilon_3-\epsilon_2)}
\right) 
\end{eqnarray}
where we use now $f(P) = P^2e^{-P_{\rm tot}/P}$ and expand it in a Taylor series around $P_{1,{\rm exp}}$ to write
\begin{eqnarray}
&&  \lim_{\epsilon_2,\epsilon_3 \to 0}\left(\frac{f(P_{1,{\rm exp}})}{\epsilon_2\epsilon_3} + 
\frac{f(P_{1,{\rm exp}})+\epsilon_2f'(P_{1,{\rm exp}})+(\epsilon_2^2/2)f''(P_{1,{\rm exp}}) + \mathcal{O}(\epsilon_2^3)}
{\epsilon_2(\epsilon_2-\epsilon_3)} +\right.\nonumber\\
&&\left. +\frac{f(P_{1,{\rm exp}})+\epsilon_3f'(P_{1,{\rm exp}})+(\epsilon_3^2/2)f''(P_{1,{\rm exp}}) + 
\mathcal{O}(\epsilon_3^3)}{\epsilon_3(\epsilon_3-\epsilon_2)} 
\right)\nonumber\\
&=& \lim_{\epsilon_2,\epsilon_3 \to 0}\left\{f(P_{1,{\rm exp}})\left[\frac{1}{\epsilon_2\epsilon_3} + 
\frac{1}{\epsilon_2(\epsilon_2-\epsilon_3)} + \frac{1}{\epsilon_3(\epsilon_3-\epsilon_2)}\right]+ \right.\nonumber\\
&&\left. + f'(P_{1,{\rm exp}})\left[\frac{\epsilon_2}{\epsilon_2(\epsilon_2-\epsilon_3)} + 
\frac{\epsilon_3}{\epsilon_3(\epsilon_3-\epsilon_2)}\right] + \right.\nonumber\\
&&\left.+ \frac{f''(P_{1,{\rm exp}})}{2}\left[\frac{\epsilon_2^2}{\epsilon_2(\epsilon_2-\epsilon_3)} + 
\frac{\epsilon_3^2}{\epsilon_3(\epsilon_3-\epsilon_2)}\right] + \mathcal{O}(\epsilon_2,\epsilon_3)\right\}\,.
\label{eq:ex3}
\end{eqnarray}
Let us now analyze the first three terms in eq. (\ref{eq:ex3}) separately.

In the first term, a simple calculation shows that the term between brackets is equal to zero. (The general argument for any value of $n$ can be seen in Section \ref{sec:application}.) 

The second term between brackets can also be shown to be zero with a simple calculation, but here we give an explicit example of the argument given in Section \ref{sec:application} by noting that it can be written as
$(1-1)/(\epsilon_2-\epsilon_3)$
where the denominator is the determinant of the $2\times 2$ Vandermonde matrix, 
$\begin{array}{|cc|}\epsilon_2 & \epsilon_3\\1 & 1\end{array}\,$, 
while the numerator is the determinant of the modified matrix resulting with one repeated line, 
$\begin{array}{|cc|}1 & 1\\1 & 1\end{array}\,$.

Lastly, we recognize that the third term between brackets has the same form as eq. (\ref{eq:lemma1-n1}) and is therefore equal to one. So we finally have the only remaining term
\begin{equation}
 \frac{f''(P_{1,{\rm exp}})}{2} 
= e^{-P_{\rm tot}/P_{1,{\rm exp}}}\left[1 + \frac{P_{\rm tot}}{P_{1,{\rm exp}}} + \frac{1}{2}\left( \frac{P_{\rm tot}}{P_{1,{\rm exp}}}\right)^2 \right]\,.
\end{equation}

\section{The ($n-1$)th-derivative of the auxiliary function $f(P)$}
\label{sec:derivative}

Here we present a proof by induction of the general expression given in eq. (\ref{eq:df}) for the $(n-1)$-derivative of the function $f(P)$ given in eq. (\ref{eq:f}). For clarity we will now index $f(P)$ by $n$ as 
\begin{equation}
f_n(P) = P^{n-1}e^{-P_{\rm tot}/P}\,.
\label{eq:fn}
\end{equation}
and the general expression for the $(n-1)$-derivative of $f_n(P)$ is reproduced below
\begin{equation}
f_n^{(n-1)}(P) = (n-1)!\, e^{-P_{\rm tot}/P}\sum_{i=0}^{n-1}\frac{(P_{\rm tot}/P)^i}{i!}\,.
\label{eq:dfn}
\end{equation}
Let us first verify that this expression is true for $n = 1$ and $n = 2$.

For $n = 1$, the trivial case, the definition (\ref{eq:fn}) gives $f_1(P) = e^{-P_{\rm tot}/P}$ and eq. (\ref{eq:dfn}) results again in $f_1(P) = e^{-P_{\rm tot}/P}$.

For $n = 2$, we have $f_2(P) = Pe^{-P_{\rm tot}/P}$ and eq.(\ref{eq:dfn}) results in $f'_2(P) = e^{-P_{\rm tot}/P}\left(1+P_{\rm tot}/P\right)$. 

Now let us assume that eq. (\ref{eq:dfn}) is true for given $n$ and $n-1$. For $n+1$, we have $f_{n+1}(P) = P^ne^{-P_{\rm tot}/P}$ and we proceed to calculate
\begin{eqnarray}
f_{n+1}^{(n)}(P) &=& \frac{d^n}{dP^n}(f_{n+1}) = \frac{d^{n-1}}{dP^{n-1}}\left(\frac{d}{dP}(f_{n+1})\right) \nonumber \\ 
&=& \frac{d^{n-1}}{dP^{n-1}}\left(nP^{n-1}e^{-P_{\rm tot}/P} + P_{\rm tot}P^{n-2}e^{-P_{\rm tot}/P}\right) \nonumber \\ 
&=& n\frac{d^{n-1}}{dP^{n-1}}(f_n(P)) + P_{\rm tot} \frac{d}{dP}\left(\frac{d^{n-2}}{dP^{n-2}}(f_{n-1}(P))\right) \nonumber\\
&=& nf_n^{(n-1)}(P) + P_{\rm tot} \frac{d}{dP}\left(f_{n-1}^{(n-2)}(P)\right)\,,
\label{eq:dfn_part}
\end{eqnarray}
where the first term is given directly by eq. (\ref{eq:dfn}) and $f_{n-1}^{(n-2)}(P)$ is given by
\begin{equation}
f_{n-1}^{(n-2)}(P) = (n-2)!\, e^{-P_{\rm tot}/P}\sum_{i=0}^{n-2}\frac{(P_{\rm tot}/P)^i}{i!}\,,
\end{equation}
for we are assuming that eq. (\ref{eq:dfn}) holds for $n$ and $n-1$. Now we calculate the first derivative of $f_{n-1}^{(n-2)}(P)$,
\begin{eqnarray}
\frac{d}{dP}\left(f_{n-1}^{(n-2)}(P)\right) &=& (n-2)!\, e^{-P_{\rm tot}/P}
\left(\sum_{i=0}^{n-2}\frac{1}{i!}\frac{P_{\rm tot}^{i+1}}{P^{i+2}}
-\sum_{i=1}^{n-2}\frac{1}{(i-1)!}\frac{P_{\rm tot}^{i}}{P^{i+1}}\right) \nonumber\\
&=& e^{-P_{\rm tot}/P}\frac{P_{\rm tot}^{n-1}}{P^n}\,,
\label{eq:dfn-2}
\end{eqnarray}
and, substituting (\ref{eq:dfn}) and (\ref{eq:dfn-2}) in  (\ref{eq:dfn_part}) we obtain
\begin{eqnarray}
f_{n+1}^{(n)}(P) &=& n!\, e^{-P_{\rm tot}/P}\sum_{i=0}^{n-1}\frac{(P_{\rm tot}/P)^i}{i!} + 
e^{-P_{\rm tot}/P}\left(\frac{P_{\rm tot}}{P}\right)^n \nonumber \\
&=& n!\, e^{-P_{\rm tot}/P}\sum_{i=0}^{n}\frac{(P_{\rm tot}/P)^i}{i!}\,,
\end{eqnarray}
completing the demonstration.


\begin{thebibliography}{1}

\bibitem{Newton} Newton, I. (1687).  \textit{Philosophae Naturalis Principia Mathematica}, London, 1687, Liber III, p. 582, Lemma V

\bibitem{LS1} Lascoux, A., Sch{\"u}tzenberger, M. P. (1982).  Polyn{\^o}mes de Schubert.  \textit{C. R. Acad. Sci. Paris. S{\'e}r. I Math.} 294: 447--450.

\bibitem{LS2} Lascoux, A., Sch{\"u}tzenberger, M. P. (1983).  Symmetry and flag manifolds.  \textit{Springer Lecture Notes in Math.} 996: 118--144.

\bibitem{LS3} Lascoux, A., Sch{\"u}tzenberger, M. P. (1987).  Symmetrizing operators on polynomial rings.  \textit{Functional Analysis and Appl.} 21: 77--78.

\bibitem{LS4} Lascoux, A., Sch{\"u}tzenberger, M. P. (1992).  D{\ ́e}compositions dans l’alg{\ ́e}bre des differences divis{\ ́e}es.  \textit{Discrete Math.} 99: 165--179.

\bibitem{Pragacz} Pragacz, P. (1996).  arXiv:alg-geom/9605014.

\bibitem{Akkouchi} Akkouchi, M. (2008).  On the convolution of exponential distributions.  \textit{J. ChungCheong Math. Soc.} 21: 501--510.

\bibitem{NumericalRecipes} Press, W. H., Teukolsky, S. A., Vetterling, W. T., Flannery,  B. P. (1992). \textit{Numerical Recipes in C: The Art of Scientific Computing, 2nd ed.} Cambridge: Cambridge University Press, pp. 90-92.

\bibitem{Schur} Schur, J. (1911). \"Uber die Darstellung der symmetrischen und der alternierenden Gruppe durch gebrochene lineare Substitutionen. \textit{J. rein angew. Math} 139: 155--250. DOI: https://doi.org/10.1515/crll.1911.139.155

\bibitem{BlackmanTukey} Blackman, R. B., Tukey, J. W. (1959). \textit{The Measurement of Power Spectra.} New York: Dover, pp. 15-25.

\bibitem{Groth75} Groth, E. J. (1975). Probability distributions related to power spectra. \textit{Astrophys. J. Suppl.} 29(286): 285--302. 


 
\end{thebibliography}
\end{document}